# WEIGHTED APPROXIMATIONS OF TAIL COPULA PROCESSES WITH APPLICATION TO TESTING THE BIVARIATE EXTREME VALUE CONDITION

By John H. J. Einmahl, Laurens de Haan and Deyuan Li[1]

*Tilburg University, Erasmus University and University of Bern*

Consider $n$ i.i.d. random vectors on $\mathbb{R}^2$, with unknown, common distribution function $F$. Under a sharpening of the extreme value condition on $F$, we derive a weighted approximation of the corresponding tail copula process. Then we construct a test to check whether the extreme value condition holds by comparing two estimators of the limiting extreme value distribution, one obtained from the tail copula process and the other obtained by first estimating the spectral measure which is then used as a building block for the limiting extreme value distribution. We derive the limiting distribution of the test statistic from the aforementioned weighted approximation. This limiting distribution contains unknown functional parameters. Therefore, we show that a version with estimated parameters converges weakly to the true limiting distribution. Based on this result, the finite sample properties of our testing procedure are investigated through a simulation study. A real data application is also presented.

**1. Introduction.**

1.1. *The extreme value model and its use.* Let $(X,Y), (X_1,Y_1), \ldots, (X_n, Y_n)$ be i.i.d. random vectors with continuous distribution function (d.f.) $F$. Suppose that there exist norming constants $a_n, c_n > 0$ and $b_n, d_n \in \mathbb{R}$ such that the sequence of d.f.'s

$$P\left(\frac{\max_{1 \leq i \leq n} X_i - b_n}{a_n} \leq x, \; \frac{\max_{1 \leq i \leq n} Y_i - d_n}{c_n} \leq y\right)$$

Received August 2004; revised September 2005.
[1]Supported in part by the Swiss National Science Foundation. Research performed mainly at Erasmus University, Rotterdam.

converges to a limit d.f., say, $G(x,y)$, with nondegenerate marginal d.f., that is,

(1.1) $$\lim_{n\to\infty} F^n(a_n x + b_n, c_n y + d_n) = G(x,y)$$

for all but countably many $x$ and $y$. Then, for a suitable choice of $a_n, b_n, c_n$ and $d_n$, there exist $\gamma_1, \gamma_2 \in \mathbb{R}$ such that

$$G(x, \infty) = \exp(-(1+\gamma_1 x)^{-1/\gamma_1}), \qquad G(\infty, y) = \exp(-(1+\gamma_2 y)^{-1/\gamma_2}).$$

The d.f. $G$ is called an extreme value d.f. and $\gamma_1, \gamma_2$ are called the (marginal) extreme value indices.

Any extreme value d.f. $G$ can be represented as

(1.2) $$G\left(\frac{x^{-\gamma_1}-1}{\gamma_1}, \frac{y^{-\gamma_2}-1}{\gamma_2}\right) \\ = \exp\left(-\int_0^{\pi/2} (x(1 \wedge \tan\theta)) \vee (y(1 \wedge \cot\theta)) \Phi(d\theta)\right),$$

with $\Phi$ the d.f. of the so-called spectral measure. There is a one-to-one correspondence between extreme value d.f.'s $G$ and finite measures with d.f. $\Phi$ that satisfy

$$\int_0^{\pi/2} (1 \wedge \tan\theta) \Phi(d\theta) = \int_0^{\pi/2} (1 \wedge \cot\theta) \Phi(d\theta) = 1,$$

via (1.2).

Alternatively, one can characterize the extreme value d.f.'s $G$ by the following: there is a measure $\Lambda$ on $[0,\infty]^2 \setminus \{(\infty,\infty)\}$ such that, with

(1.3) $$l(x,y) := -\log G\left(\frac{x^{-\gamma_1}-1}{\gamma_1}, \frac{y^{-\gamma_2}-1}{\gamma_2}\right),$$

we have

(1.4) 
1. $l(x,y) = \Lambda(\{(u,v) \in [0,\infty]^2 : u \leq x \text{ or } v \leq y\})$,
2. $l(ax, ay) = al(x,y) \qquad$ for $a, x, y > 0$.

More generally, we have for any $a > 0$ and any Borel set $A \subset [0,\infty]^2 \setminus \{(\infty,\infty)\}$

(1.5) $$\Lambda(aA) = a\Lambda(A),$$

with $aA := \{(ax, ay) : (x,y) \in A\}$. Also, (1.1) implies

(1.6) $$\lim_{t \downarrow 0} t^{-1} P((1-F_1(X), 1-F_2(Y)) \in tA) = \Lambda(A)$$

for any Borel set $A$, provided $\Lambda(\partial A) = 0$, where $F_1(x) := F(x,\infty)$ and $F_2(y) := F(\infty, y)$. See, for example, [2].

The bivariate extreme value framework is the appropriate one when one wants to estimate the probability of an *extreme set*, that is, a set outside the range of even the largest observations. Take $a > 0$ small. Since by (1.5) and (1.6), for small $t$,

$$P((1 - F_1(X), 1 - F_2(Y)) \in taA) \approx aP((1 - F_1(X), 1 - F_2(Y)) \in tA),$$

we can estimate the probability of $tA$—outside the range of the observations—asymptotically by estimating the probability of the pulled back set $taA$ using the empirical measure. See [3]. Condition (1.1) is fulfilled for many standard distributions but not for all distributions. Hence, before using this framework to estimate probabilities of extreme sets, it is important to check whether (1.1) is a reasonable assumption for the data set at hand. And one wants to do this beforehand, without specifying the exact structure of the limiting distribution.

1.2. *Estimation of model parameters.* Now, in order to develop a test, let us consider the following. Relation (1.1) implies [cf. (1.6)]

$$\lim_{t \downarrow 0} t^{-1} P((1 - F_1(X)) \wedge (1 - F_2(Y)) \leq t,$$

(1.7)

$$1 - F_2(Y) \leq (1 - F_1(X)) \tan \theta) = \Phi(\theta)$$

for continuity points $\theta \in (0, \pi/2]$ of $\Phi$. Also, for $(x, y) \in [0, \infty)^2$,

(1.8) $$\lim_{t \downarrow 0} t^{-1} P(1 - F_1(X) \leq tx \text{ or } 1 - F_2(Y) \leq ty) = l(x, y).$$

A nonparametric estimator for $\Phi$, suggested by the limit relation (1.7) is [8]

(1.9) $$\hat{\Phi}(\theta) := \frac{1}{k} \sum_{i=1}^{n} I_{\{R_i^X \vee R_i^Y \geq n+1-k,\ n+1-R_i^Y \leq (n+1-R_i^X) \tan \theta\}},$$

where $R_i^X$ is the rank of $X_i$ among $X_1, X_2, \ldots, X_n$ and $R_i^Y$ is the rank of $Y_i$ among $Y_1, Y_2, \ldots, Y_n$, where $k = k(n)$ is an intermediate sequence of integers, that is, $k \to \infty$, $k/n \to 0$, as $n \to \infty$. Similarly, a nonparametric estimator for $l$, suggested by the limit relation (1.8), is ([9]; see also [6])

$$\hat{l}_2(x, y) := \frac{1}{k} \sum_{i=1}^{n} I_{\{X_i > X_{n+1-\lceil kx \rceil:n} \text{ or } Y_i > Y_{n+1-\lceil ky \rceil:n}\}}$$

(1.10)

$$= \frac{1}{k} \sum_{i=1}^{n} I_{\{R_i^X > n+1-kx \text{ or } R_i^Y > n+1-ky\}},$$

where $X_{1:n} \leq \cdots \leq X_{n:n}$ are the order statistics of the $X_i$, $i = 1, 2, \ldots, n$ (similarly for the $Y_i$), with $\lceil z \rceil$ the smallest integer $\geq z$.

Another way of estimating $l$ is via (1.2), (1.3) and (1.9):

$$(1.11) \qquad \hat{l}_1(x,y) := \int_0^{\pi/2} (x(1 \wedge \tan\theta)) \vee (y(1 \wedge \cot\theta)) \hat{\Phi}(d\theta).$$

1.3. *The test.* A promising approach to testing whether the null hypothesis (1.8) holds seems to be to see if the two estimators $\hat{l}_1$ and $\hat{l}_2$ for $l$, that have a different background, are not too different. The estimator $\hat{l}_2$ is a natural one mimicking more or less the tail of the distribution itself. But this estimator does not necessarily satisfy condition 2 of (1.4). On the other hand, $\hat{l}_1$ does satisfy condition 2 of (1.4), but the estimator itself is of a somewhat more complicated nature.

The proposed test statistic is of the Anderson–Darling type:

$$(1.12) \qquad L_n := \iint_{0 < x,y \leq 1} (\hat{l}_1(x,y) - \hat{l}_2(x,y))^2 (x \vee y)^{-\beta} \, dx \, dy$$

for certain $\beta \geq 0$. The test statistic is similar to those used for testing a parametric null hypothesis (like testing for normality), where the empirical distribution function is compared with the true distribution function with estimated parameters. Here, however, the estimated parameter $\Phi$ is a function (and we only deal with the tail of the distribution). Also, note that our methods allow us to deal with test statistics other than $L_n$ as well.

It is not difficult to see that if relation (1.8) is true, the statistic $L_n$ tends to zero in probability as $n \to \infty$. We shall establish the asymptotic distribution of $kL_n$ as $n \to \infty$ under (1.1) and some extra conditions stemming from [9] and [8], thus providing a basis for applying a test. The hypothesis to be tested is (1.8). For the asymptotic normality of the test statistic $kL_n$, under $H_0$, extra conditions are needed. See Remark 2.2 below.

Note that this test checks whether the dependence structure satisfies (1.8) and not if the marginals $F_1, F_2$ are of the right type. It is only based on the relative positions (ranks) of the data and completely independent of the marginal distributions for which tests have been developed already in [5] and [4].

1.4. *Use of test.* As mentioned before, the test statistic $kL_n$ can be used for a preliminary test of the extreme value model (1.1) before one uses the model in applications. Note that the test statistic $kL_n$ is based on observations for which at least one component exceeds a certain threshold. Since the estimators depend on this threshold, one can plot $kL_n$ as a function of $k$. This plot can be used as an exploratory tool for determining from which threshold on the two estimators $\hat{l}_1$ and $\hat{l}_2$ are close to each other, suggesting that the approximations (1.7) and (1.8) can be trusted, and, hence, yields a heuristic procedure for determining $k$. So this is a second use of the test statistic $kL_n$. See also [14], Section 5.4, and [1].

1.5. *Outline of the paper.* The weak convergence of $kL_n$ is stated in Theorem 2.3. For the proof of this theorem, the known asymptotic normality result for $\hat{\Phi}$ [8] is sufficient but not the known one for $\hat{l}_2$ [9]. Hence, as a preliminary but important result, we first develop a Gaussian approximation for the weighted tail copula process on $(0,1]^2$,

$$\sqrt{k}(\hat{l}_2(x,y) - l(x,y))/(x \vee y)^\eta, \qquad 0 \leq \eta < 1/2,$$

thus extending significantly the result of [9], where $\eta = 0$. This result is stated in Theorem 2.2. The proofs are given in Section 3.

The limiting random variable in Theorem 2.3 is determined as an integral of a combination of Gaussian processes. They are parametrized by functions which can be estimated consistently. In Section 4 it is proved that the probability distribution of the limiting random variable with these functions estimated converges to the distribution of the limiting random variable with these functions equal to the actual ones, which makes the procedure applicable in practice. In Section 5 simulation results and an application to real data are reported.

**2. Main results.** Before stating the main results, we need to introduce some notation.

Let $W_\Lambda$ be a Wiener process indexed by the Borel sets in $[0,\infty]^2 \setminus \{(\infty,\infty)\}$, depending on the parameter $\Lambda$ from (1.4) in the following way: $W_\Lambda$ is a Gaussian process and for Borel sets $C$ and $\tilde{C}$,

$$EW_\Lambda(C) = 0 \quad \text{and} \quad EW_\Lambda(C)W_\Lambda(\tilde{C}) = \Lambda(C \cap \tilde{C}).$$

Define the sets $C_\theta$ by

$$C_\theta = \{(x,y) \in [0,\infty]^2 : x \wedge y \leq 1, y \leq x \tan\theta\}, \qquad \theta \in [0, \pi/2].$$

Assume that the measure $\Lambda$ has a density $\lambda$. The process $Z$ on $[0, \pi/2]$ is defined by

$$
\begin{aligned}
Z(\theta) &= \int_0^{1 \vee 1/(\tan\theta)} \lambda(x, x\tan\theta)(W_1(x)\tan\theta - W_2(x\tan\theta))\,dx \\
&\quad - W_2(1) \int_{1 \vee 1/(\tan\theta)}^\infty \lambda(x,1)\,dx \\
&\quad - I_{(\pi/4,\pi/2]}(\theta) W_1(1) \int_1^{\tan\theta} \lambda(1,y)\,dy, \qquad \theta \in [0, \pi/2), \\
Z\left(\frac{\pi}{2}\right) &= -W_2(1) \int_1^\infty \lambda(x,1)\,dx - W_1(1) \int_1^\infty \lambda(1,y)\,dy,
\end{aligned}
$$

(2.1)

where $W_1, W_2$ are the marginal processes defined by

$$W_1(x) = W_\Lambda([0,x] \times [0,\infty]) \quad \text{and} \quad W_2(y) = W_\Lambda([0,\infty] \times [0,y]).$$

Clearly, both processes are standard Wiener processes. Define, for $x, y > 0$,

(2.2) $\qquad R(x,y) = \Lambda([0,x] \times [0,y]),$

(2.3) $\qquad W_R(x,y) = W_\Lambda([0,x] \times [0,y]),$

(2.4) $\qquad R_1(x,y) = \partial R(x,y)/\partial x, \qquad R_2(x,y) = \partial R(x,y)/\partial y.$

For convenient presentation and convenient application, the next two theorems are presented in an approximation setting (with all the processes involved defined on *one* probability space), via the Skorohod construction. So in these theorems, $\hat{l}_1$, $\hat{l}_2$ and the limiting processes $A$ and $B$ (defined below) are only equal in distribution to the original ones, but we do not add the usual tildes to the notation.

THEOREM 2.1. *Assume that condition (1.6) and Conditions 1 and 2 of [8] hold, and that $\Lambda$ has a continuous density $\lambda$ on $[0,\infty)^2 \setminus \{(0,0)\}$. Then*

$$\sup_{0 < x,y \leq 1} \frac{|\sqrt{k}(\hat{l}_1(x,y) - l(x,y)) - A(x,y)|}{x \vee y} \xrightarrow{P} 0$$

*as $n \to \infty$, where*

$$A(x,y) := \begin{cases} x\left(W_\Lambda(C_{\pi/2}) + Z\left(\frac{\pi}{2}\right)\right) \\ \quad + y \int_{\pi/4}^{\arctan y/x} \frac{1}{\sin^2 \theta} (W_\Lambda(C_\theta) + Z(\theta))\, d\theta, & \text{if } y \geq x, \\ x\left(W_\Lambda(C_{\pi/2}) + Z\left(\frac{\pi}{2}\right)\right) \\ \quad - x \int_{\arctan y/x}^{\pi/4} \frac{1}{\cos^2 \theta} (W_\Lambda(C_\theta) + Z(\theta))\, d\theta, & \text{if } y < x. \end{cases}$$

Write $U_i = 1 - F_1(X_i), V_i = 1 - F_2(Y_i), i = 1, 2, \ldots, n$. Let $C(x,y)$ be the distribution function of $(U_i, V_i)$. By (1.6) and (2.2), we have

$$R(x,y) = \lim_{t \downarrow 0} t^{-1} C(tx, ty).$$

We assume, as in [9], that, for some $\alpha > 0$,

(2.5) $\qquad t^{-1} C(tx, ty) - R(x,y) = O(t^\alpha) \qquad \text{as } t \downarrow 0,$

uniformly for $x \vee y \leq 1$, $x, y \geq 0$.

THEOREM 2.2. *Assume that conditions (1.6) and (2.5) hold and that $k = o(n^{2\alpha/(1+2\alpha)})$. If $R_1$ and $R_2$ are continuous, then we have, for $0 \leq \eta < 1/2$,*

$$\sup_{0 < x,y \leq 1} \frac{|\sqrt{k}(\hat{l}_2(x,y) - l(x,y)) + B(x,y)|}{(x \vee y)^\eta} \xrightarrow{P} 0$$

as $n \to \infty$, where $B(x,y) := W_R(x,y) - R_1(x,y)W_1(x) - R_2(x,y)W_2(y)$.

THEOREM 2.3. *Assume the conditions of Theorems* 2.1 *and* 2.2 *hold. Then for each $0 \le \beta < 3$,*

(2.6)
$$\iint_{0<x,y\le 1} \frac{k(\hat{l}_1(x,y) - \hat{l}_2(x,y))^2}{(x \vee y)^\beta} \, dx \, dy$$
$$\xrightarrow{d} \iint_{0<x,y\le 1} \frac{(A(x,y) + B(x,y))^2}{(x \vee y)^\beta} \, dx \, dy$$

*as $n \to \infty$, and the limit is finite almost surely.*

REMARK 2.1. The case $\beta = 0$ is similar to the Cramér–von Mises test. Note that for $\beta < 2$, Theorem 2.3 easily follows from an unweighted approximation in Theorems 2.1 and 2.2. Therefore, the case $\beta = 2$ (and not, as usual, $\beta = 1$) is similar to the Anderson–Darling test.

REMARK 2.2. Conditions 1 and 2 of [8] are rather involved. They require that the convergence in (1.6) hold uniformly over certain classes of sets. Moreover, they put an extra restriction on the growth of the sequence $k(n)$, related to the rate of that convergence. The assumption that $\Lambda$ has a density $\lambda$ excludes, for example, asymptotic independence, that is, $l(x,y) = x + y$, for all $x, y \ge 0$. Condition (2.5) is rather weak, but there are distributions for which (1.1) holds with asymptotic dependence, but where the rate of convergence is slower than $t^\alpha$ for any $\alpha > 0$.

REMARK 2.3. The random variable on the right in Theorem 2.3 has a continuous distribution function. This follows from a property of Gaussian measures on Banach spaces: the measure of a closed ball is a continuous function of its radius; see, for example, [13], Chapter 4, Theorem 1.2.

REMARK 2.4. Since $x \vee y \le l(x,y) \le x + y \le 2(x \vee y)$, (2.6) remains true with $x \vee y$ replaced with $l(x,y)$ or $x + y$, but, when choosing $l(x,y)$, the left-hand side of (2.6) is not a statistic and $l$ has to be estimated.

**3. Proofs.** Before proving Theorem 2.1, we first present two lemmas and a proposition.

LEMMA 3.1.
$$l(x,y) = \begin{cases} x\Phi\left(\dfrac{\pi}{2}\right) + y \displaystyle\int_{\pi/4}^{\arctan y/x} \dfrac{1}{\sin^2 \theta} \Phi(\theta) \, d\theta, & \text{if } y \ge x, \\ x\Phi\left(\dfrac{\pi}{2}\right) - x \displaystyle\int_{\arctan y/x}^{\pi/4} \dfrac{1}{\cos^2 \theta} \Phi(\theta) \, d\theta, & \text{if } y < x. \end{cases}$$

PROOF. Since
$$l(x,y) = \int_0^{\pi/2} (x(1 \wedge \tan\theta)) \vee (y(1 \wedge \cot\theta)) \Phi(d\theta)$$
$$= \int_0^{\pi/4} (x\tan\theta) \vee y \Phi(d\theta) + \int_{\pi/4}^{\pi/2} x \vee (y\cot\theta) \Phi(d\theta)$$

and $x\tan\theta > y \Leftrightarrow x > y\cot\theta \Leftrightarrow \theta > \arctan\frac{y}{x}$, we have

$$l(x,y) = \int_0^{\pi/4 \wedge \arctan y/x} y\Phi(d\theta) + \int_{\pi/4 \wedge \arctan y/x}^{\pi/4} x\tan\theta \Phi(d\theta)$$
$$+ \int_{\pi/4}^{\pi/4 \vee \arctan y/x} y\cot\theta \Phi(d\theta) + \int_{\pi/4 \vee \arctan y/x}^{\pi/2} x\Phi(d\theta)$$
$$= \begin{cases} \int_0^{\pi/4} y\Phi(d\theta) + \int_{\pi/4}^{\arctan y/x} y\cot\theta \Phi(d\theta) + \int_{\arctan y/x}^{\pi/2} x\Phi(d\theta), \\ \qquad\qquad\qquad\qquad\qquad\qquad\qquad\qquad\qquad \text{if } y \geq x, \\ \int_0^{\arctan y/x} y\Phi(d\theta) + \int_{\arctan y/x}^{\pi/4} x\tan\theta \Phi(d\theta) + \int_{\pi/4}^{\pi/2} x\Phi(d\theta), \\ \qquad\qquad\qquad\qquad\qquad\qquad\qquad\qquad\qquad \text{if } y < x. \end{cases}$$

In the case $y \geq x$, via integration by parts, one has
$$l(x,y) = y\Phi\left(\frac{\pi}{4}\right) - y\Phi(0) + y\cot\left(\arctan\frac{y}{x}\right)\Phi\left(\arctan\frac{y}{x}\right) - y\cot\frac{\pi}{4}\Phi\left(\frac{\pi}{4}\right)$$
$$- y\int_{\pi/4}^{\arctan y/x} \Phi(\theta)\left(-\frac{1}{\sin^2\theta}\right) d\theta + x\Phi\left(\frac{\pi}{2}\right) - x\Phi\left(\arctan\frac{y}{x}\right)$$
$$= x\Phi\left(\frac{\pi}{2}\right) + y\int_{\pi/4}^{\arctan y/x} \frac{1}{\sin^2\theta}\Phi(\theta) d\theta.$$

In the case $y < x$, via integration by parts again, one has
$$l(x,y) = y\Phi\left(\arctan\frac{y}{x}\right) - y\Phi(0) + x\tan\frac{\pi}{4}\Phi\left(\frac{\pi}{4}\right)$$
$$- x\tan\left(\arctan\frac{y}{x}\right)\Phi\left(\arctan\frac{y}{x}\right)$$
$$- x\int_{\arctan y/x}^{\pi/4} \Phi(\theta)\frac{1}{\cos^2\theta} d\theta + x\Phi\left(\frac{\pi}{2}\right) - x\Phi\left(\frac{\pi}{4}\right)$$
$$= x\Phi\left(\frac{\pi}{2}\right) - x\int_{\arctan y/x}^{\pi/4} \frac{1}{\cos^2\theta}\Phi(\theta) d\theta. \qquad \square$$

Write
$$(3.1) \quad R_n(x,y) = \frac{n}{k}C\left(\frac{kx}{n}, \frac{ky}{n}\right), \qquad T_n(x,y) = \frac{1}{k}\sum_{i=1}^n I_{\{U_i < kx/n, V_i < ky/n\}},$$

(3.2) $\quad v_n(x,y) = \sqrt{k}(T_n(x,y) - R_n(x,y)), \qquad v_{n,\eta}(x,y) = \dfrac{v_n(x,y)}{(x \vee y)^\eta}$

and

(3.3)
$$v_{n,\eta,1}(x) = \frac{v_n(x,\infty)}{x^\eta},$$

$$v_{n,\eta,2}(y) = \frac{v_n(\infty,y)}{y^\eta}, \qquad v_{n,j} = v_{n,0,j}, j=1,2.$$

PROPOSITION 3.1. *Let $T > 0$. For $0 \leq \eta < 1/2$,*

$$(v_{n,\eta}(x,y), x, y \in (0,T], v_{n,\eta,1}(x), x \in (0,T], v_{n,\eta,2}(y), y \in (0,T])$$

*converges in distribution to*

$$\left( \frac{W_R(x,y)}{(x \vee y)^\eta}, x, y \in (0,T], \frac{W_1(x)}{x^\eta}, x \in (0,T], \frac{W_2(y)}{y^\eta}, y \in (0,T] \right)$$

*as $n \to \infty$.*

PROOF. Define

$$Z_{n,i} = \frac{1}{\sqrt{k}} \delta_{((n/k)U_i, (n/k)V_i)},$$

and for all $0 < x, y \leq T$, define the functions

$$f_{x,y} = I_{[0,x) \times [0,y)}/(x \vee y)^\eta,$$

$$f_x^{(1)} = I_{[0,x) \times [0,\infty]}/x^\eta, \qquad f_y^{(2)} = I_{[0,\infty] \times [0,y)}/y^\eta.$$

All these $f$'s form the class $\mathcal{F}$. We equip $\mathcal{F}$ with the semi-metric $d$ defined by

$$d(f_{x,y}, f_{u,v}) = \sqrt{E\left( \frac{W_R(x,y)}{(x \vee y)^\eta} - \frac{W_R(u,v)}{(u \vee v)^\eta} \right)^2},$$

$$d(f_{x,y}, f_u^{(1)}) = \sqrt{E\left( \frac{W_R(x,y)}{(x \vee y)^\eta} - \frac{W_1(u)}{u^\eta} \right)^2},$$

and so on.

For any $\varepsilon > 0$, the bracketing number $N_{[]}(\varepsilon, \mathcal{F})$ is the minimal number of sets $N_\varepsilon$ in a partition $\mathcal{F} = \bigcup_{j=1}^{N_\varepsilon} \mathcal{F}_{\varepsilon j}$ of the index set into sets $\mathcal{F}_{\varepsilon j}$ such that, for every partitioning set $\mathcal{F}_{\varepsilon j}$,

(3.4) $$\sum_{i=1}^n E^* \sup_{f,g \in \mathcal{F}_{\varepsilon j}} |Z_{n,i}(f) - Z_{n,i}(g)|^2 \leq \varepsilon^2,$$

where, as usual, $Z_{n,i}(f) = \int f \, dZ_{n,i}$ and $E^*$ means taking the outer integral when computing the expectation.

We will use Theorem 2.11.9 in [17]: For each $n$, let $Z_{n,1}, Z_{n,2}, \ldots, Z_{n,n}$ be independent stochastic processes with finite second moments indexed by a totally bounded semimetric space $(\mathcal{F}, d)$. Suppose

$$\sum_{i=1}^{n} E^* \|Z_{n,i}\|_{\mathcal{F}} \mathbb{1}_{\{\|Z_{n,i}\|_{\mathcal{F}} > \lambda\}} \to 0 \qquad \text{for every } \lambda > 0,$$

where $\|Z_{n,i}\|_{\mathcal{F}} = \sup_{f \in \mathcal{F}} |Z_{n,i}(f)|$, and

$$\int_0^{\delta_n} \sqrt{\log N_{[]}(\varepsilon, \mathcal{F})} \, d\varepsilon \to 0 \qquad \text{for every } \delta_n \downarrow 0.$$

Then the sequence $\sum_{i=1}^{n}(Z_{n,i} - EZ_{n,i})$ is asymptotically tight in $\ell^\infty(\mathcal{F})$ and converges weakly, provided the finite-dimensional distributions converge weakly.

We briefly sketch the total boundedness of $(\mathcal{F}, d)$. We only consider the subclass $\mathcal{F}_2$ of $\mathcal{F}$ consisting of the bivariate $f_{x,y}$'s; moreover, we restrict ourselves to the case $x \geq y$, $u \geq v$ and $x \geq u$, $y \geq v$. For any $\delta > 0$, assuming $|x - u| \leq \delta$ and $|y - v| \leq \delta$, one has

$$\begin{aligned} d^2(f_{x,y}, f_{u,v}) &= E\left(\frac{W_R(x,y)}{(x \vee y)^\eta} - \frac{W_R(u,v)}{(u \vee v)^\eta}\right)^2 \\ &= E\left(\frac{u^\eta W_R(x,y) - x^\eta W_R(u,v)}{(xu)^\eta}\right)^2 \\ &= \frac{u^{2\eta} R(x,y) - 2x^\eta u^\eta R(u,v) + x^{2\eta} R(u,v)}{(xu)^{2\eta}}. \end{aligned}$$

If $u \leq \delta$, then

$$\begin{aligned} d^2(f_{x,y}, f_{u,v}) &\leq \frac{R(x,y)}{x^{2\eta}} + \frac{2R(u,v)}{u^{2\eta}} + \frac{R(u,v)}{u^{2\eta}} \\ &\leq x^{1-2\eta} + 3u^{1-2\eta} \leq (2\delta)^{1-2\eta} + 3\delta^{1-2\eta} \leq 5\delta^{1-2\eta}. \end{aligned}$$

If $u > \delta$, then, since

$$R(x,y) \leq R(u,v) + \Lambda([u,x] \times [0,\infty]) + \Lambda([0,\infty] \times [v,y]) \leq R(u,v) + 2\delta,$$

we have

$$\begin{aligned} d^2(f_{x,y}, f_{u,v}) &\leq \frac{R(u,v)(u^\eta - x^\eta)^2}{(xu)^{2\eta}} + \frac{2\delta u^{2\eta}}{(xu)^{2\eta}} \\ &\leq u^{1-4\eta}(u^\eta - x^\eta)^2 + 2\delta^{1-2\eta} \\ &\leq u^{1-4\eta} x^{2\eta-2}(x-u)^2 + 2\delta^{1-2\eta} \\ &\leq u^{-1-2\eta}(x-u)^2 + 2\delta^{1-2\eta} \leq 3\delta^{1-2\eta}. \end{aligned}$$

So, since $1 - 2\eta > 0$, we see that, for every $\varepsilon > 0$, we can find a $\delta > 0$ such that, for $|x - u| \leq \delta$ and $|y - v| \leq \delta$, $d^2(f_{x,y}, f_{u,v}) < \varepsilon$. Hence, since $[0, T]^2$ is totally bounded with respect to the Euclidean metric, we obtain the total boundedness of $(\mathcal{F}, d)$.

Observe that
$$Z_{n,i}(f_{x,y}) = \frac{1}{\sqrt{k}} I_{\{U_i < (k/n)x, V_i < (k/n)y\}}/(x \vee y)^\eta,$$

$$\sum_{i=1}^n (Z_{n,i} - EZ_{n,i})(f_{x,y}) = v_{n,\eta}(x,y)$$

and similarly for the marginal processes. First we have to show that, for every $\lambda > 0$,

(3.5) $$\sum_{i=1}^n E\|Z_{n,i}\|_\mathcal{F} I_{\{\|Z_{n,i}\|_\mathcal{F} > \lambda\}} \to 0$$

as $n \to \infty$. Again, we will restrict ourselves to the subclass $\mathcal{F}_2$. For the univariate $f_x^{(1)}$'s and $f_y^{(2)}$'s, it can be shown in a similar but easier way.

Note that
$$\sup_{f_{x,y} \in \mathcal{F}_2} \frac{1}{\sqrt{k}} I_{\{U_i < (k/n)x, V_i < (k/n)y\}}/(x \vee y)^\eta \leq \frac{1}{\sqrt{k}} \frac{1}{((n/k)(U_i \vee V_i))^\eta},$$

so for each $\lambda > 0$,
$$\sum_{i=1}^n E\|Z_{n,i}\|_{\mathcal{F}_2} I_{\{\|Z_{n,i}\|_{\mathcal{F}_2} > \lambda\}}$$

$$\leq \frac{n}{\sqrt{k}} E \frac{1}{((n/k)(U_1 \vee V_1))^\eta} I_{\{(n/k)(U_1 \vee V_1) < (\sqrt{k}\lambda)^{-1/\eta}\}}$$

$$= \frac{n}{\sqrt{k}} \int_0^{(\sqrt{k}\lambda)^{-1/\eta}} x^{-\eta} \, dC\left(\frac{k}{n}x, \frac{k}{n}x\right)$$

$$= \frac{n}{\sqrt{k}} \left(\sqrt{k}\lambda C\left(\frac{k}{n}(\sqrt{k}\lambda)^{-1/\eta}, \frac{k}{n}(\sqrt{k}\lambda)^{-1/\eta}\right)\right.$$

$$\left. + \eta \int_0^{(\sqrt{k}\lambda)^{-1/\eta}} C\left(\frac{k}{n}x, \frac{k}{n}x\right) x^{-\eta-1} \, dx\right)$$

$$\leq \frac{n}{\sqrt{k}} \left(\sqrt{k}\lambda \frac{k}{n}(\sqrt{k}\lambda)^{-1/\eta} + \eta \int_0^{(\sqrt{k}\lambda)^{-1/\eta}} \frac{k}{n} x^{-\eta} \, dx\right)$$

$$= \lambda^{1-1/\eta} k^{1-1/(2\eta)} + \sqrt{k} \frac{\eta}{1-\eta} (\sqrt{k}\lambda)^{1-1/\eta}$$

$$= \frac{1}{1-\eta} \lambda^{1-1/\eta} k^{1-1/(2\eta)} \to 0 \qquad (\eta < 1/2).$$

Next, we want to prove

(3.6) $$\int_0^{\delta_n} \sqrt{\log N_{[]}(\varepsilon, \mathcal{F})}\, d\varepsilon \to 0$$

for every $\delta_n \downarrow 0$. For notational convenience, we choose $T = 1$; for general $T > 0$, the proof goes the same. Let $\varepsilon > 0$ be small, define $a = \varepsilon^{3/(1-2\eta)}$ and $\theta = 1 - \varepsilon^3$. We again consider only $\mathcal{F}_2$; the univariate $f$'s are easier to handle. Define

$$\mathcal{F}(a) = \{f_{x,y} \in \mathcal{F}_2 : x \wedge y \leq a\},$$
$$\mathcal{F}(l,m) = \{f_{x,y} \in \mathcal{F}_2 : \theta^{l+1} \leq x \leq \theta^l, \theta^{m+1} \leq y \leq \theta^m\}.$$

Then

$$\mathcal{F}_2 = \mathcal{F}(a) \cup \left(\bigcup_{m=0}^{[\log a/\log \theta]} \bigcup_{l=0}^{[\log a/\log \theta]} \mathcal{F}(l,m)\right).$$

First check (3.4) for $\mathcal{F}(a)$:

$$\sum_{i=1}^n E \sup_{f,g \in \mathcal{F}(a)} (Z_{n,i}(f) - Z_{n,i}(g))^2$$
$$= nE \sup_{f,g \in \mathcal{F}(a)} (Z_{n,1}(f) - Z_{n,1}(g))^2$$
$$\leq 4nE \sup_{f \in \mathcal{F}(a)} Z_{n,1}^2(f)$$
$$= \frac{4n}{k} E \sup_{\substack{x,y>0 \\ x \wedge y \leq a}} I_{\{U_1 < kx/n, V_1 < ky/n\}}/(x \vee y)^{2\eta}$$
$$\leq \frac{4n}{k} E \left(\frac{n}{k} U_1\right)^{-2\eta} I_{\{(n/k)U_1 < a\}}$$
$$= \frac{4n}{k} \int_0^{ak/n} \left(\frac{n}{k} x\right)^{-2\eta} dx$$
$$= \frac{4}{1-2\eta} a^{1-2\eta} \leq \varepsilon^2.$$

Now we consider (3.4) for the $\mathcal{F}(l,m)$; w.l.o.g. we take $l \leq m$:

$$\sum_{i=1}^n E \sup_{f,g \in \mathcal{F}(l,m)} (Z_{n,i}(f) - Z_{n,i}(g))^2$$
$$\leq nE\left(\sup_{f \in \mathcal{F}(l,m)} Z_{n,1}(f) - \inf_{f \in \mathcal{F}(l,m)} Z_{n,1}(f)\right)^2$$

$$\leq \frac{n}{k} E(I_{\{U_1<(k/n)\theta^l, V_1<(k/n)\theta^m\}}/(\theta^{l+1} \vee \theta^{m+1})^\eta$$

$$- I_{\{U_1<(k/n)\theta^{l+1}, V_1<(k/n)\theta^{m+1}\}}/(\theta^l \vee \theta^m)^\eta)^2$$

$$= \frac{n}{k} E\bigg(I_{\{U_1<(k/n)\theta^l, V_1<(k/n)\theta^m\}}\bigg(\frac{1}{\theta^{\eta(l+1)}} - \frac{1}{\theta^{\eta l}}\bigg)$$

$$+ (I_{\{U_1<(k/n)\theta^l, V_1<(k/n)\theta^m\}} - I_{\{U_1<(k/n)\theta^{l+1}, V_1<(k/n)\theta^{m+1}\}})\frac{1}{\theta^{\eta l}}\bigg)^2$$

$$\leq \frac{2n}{k}\bigg(C\bigg(\frac{k}{n}\theta^l, \frac{k}{n}\theta^m\bigg)\frac{1}{\theta^{2\eta l}}\bigg(\frac{1}{\theta^\eta} - 1\bigg)^2$$

$$+ \bigg[C\bigg(\frac{k}{n}\theta^l, \frac{k}{n}\theta^m\bigg) - C\bigg(\frac{k}{n}\theta^{l+1}, \frac{k}{n}\theta^{m+1}\bigg)\bigg]\frac{1}{\theta^{2\eta l}}\bigg)$$

$$\leq \frac{2n}{k}\bigg(\frac{k}{n}\frac{\theta^l}{\theta^{2\eta l}}\bigg(\frac{1}{\theta^\eta} - 1\bigg)^2 + \frac{2k}{n}\frac{\theta^l}{\theta^{2\eta l}}(1-\theta)\bigg)$$

$$\leq 2\bigg(\frac{1}{\theta^{1/2}} - 1\bigg)^2 + 4(1-\theta) \leq \varepsilon^6 + 4\varepsilon^3 \leq \varepsilon^2.$$

It is easy to see that the number of elements of the "partition" of $\mathcal{F}_2$ is bounded by $\varepsilon^{-7}$, which yields (3.6). Hence, we have proved the asymptotic tightness condition.

It remains to prove that the finite-dimensional distributions of our process converge weakly. This follows from the fact that multivariate weak convergence follows from weak convergence of linear combinations of the components and the univariate Lindeberg–Feller central limit theorem. It is easily seen that the Lindeberg condition is satisfied for these linear combinations since the elements of $\mathcal{F}$ are weighted indicators and, hence, bounded. □

LEMMA 3.2. *For $0 \leq \eta < 1/2$,*

$$P\bigg(\sup_{\substack{x \vee y \leq \varepsilon \\ x,y>0}} \frac{|W_R(x,y)|}{(x \vee y)^\eta} \geq \lambda\bigg) \leq 16 \sum_{m=0}^\infty \exp\bigg(-\frac{\lambda^2}{2^{1+2\eta}}\frac{2^{m(1-2\eta)}}{\varepsilon^{1-2\eta}}\bigg).$$

PROOF. For $m = 0, 1, 2, \ldots$, define

$$\mathcal{A}_m = \bigg\{(x,y): \frac{\varepsilon}{2^{m+1}} \leq x \leq \frac{\varepsilon}{2^m}, \frac{\varepsilon}{2^{m+1}} \leq y \leq \varepsilon\bigg\}.$$

Then, with $Z$ a standard normal random variable,

$$P\bigg(\sup_{\substack{x \vee y \leq \varepsilon \\ 0 < x \leq y}} \frac{|W_R(x,y)|}{(x \vee y)^\eta} \geq \lambda\bigg) = P\bigg(\sup_{\substack{x \vee y \leq \varepsilon \\ 0 < x \leq y}} \frac{|W_R(x,y)|}{y^\eta} \geq \lambda\bigg)$$

$$\leq P\bigg(\sup_{m\in\{0,1,2,\ldots\}}\sup_{(x,y)\in\mathcal{A}_m}\frac{|W_R(x,y)|}{y^\eta}\geq\lambda\bigg)$$

$$\leq \sum_{m=0}^{\infty} P\bigg(\sup_{(x,y)\in\mathcal{A}_m}|W_R(x,y)|\geq\lambda\Big(\frac{\varepsilon}{2^{m+1}}\Big)^\eta\bigg)$$

$$\leq 4\sum_{m=0}^{\infty} P\bigg(\Big|W_R\Big(\frac{\varepsilon}{2^m},\varepsilon\Big)\Big|\geq\lambda\Big(\frac{\varepsilon}{2^{m+1}}\Big)^\eta\bigg)$$

$$\leq 4\sum_{m=0}^{\infty} P\bigg(|Z|\geq\frac{\lambda}{2^\eta}\Big(\frac{2^m}{\varepsilon}\Big)^{1/2-\eta}\bigg)$$

$$\leq 8\sum_{m=0}^{\infty}\exp\bigg(-\frac{\lambda^2}{2^{1+2\eta}}\frac{2^{m(1-2\eta)}}{\varepsilon^{1-2\eta}}\bigg),$$

where the third inequality follows, for instance, from an adaptation of Lemma 1.2 in [12] and the last inequality from Mill's ratio. A symmetry argument completes the proof. □

By Theorem 2 in [8] and Proposition 3.1 (and their proofs), it follows that

$$(\sqrt{k}(\hat{\Phi}(\theta)-\Phi(\theta)),v_{n,\eta}(x,y),v_{n,\eta,1}(u),v_{n,\eta,2}(v))$$
$$\xrightarrow{d}\bigg(W_\Lambda(C_\theta)+Z(\theta),\frac{W_R(x,y)}{(x\vee y)^\eta},\frac{W_1(u)}{u^\eta},\frac{W_2(v)}{v^\eta}\bigg)$$

on $D[0,\pi/2]\times D[0,T]^2\times D[0,T]\times D[0,T]$. By the Skorohod construction, there exists now a probability space carrying $\hat{\Phi}^*$, $v_n^*$, $v_{n,1}^*$, $v_{n,2}^*$, $W_\Lambda^*(C.)$, $Z^*$, $W_R^*$, $W_1^*$ and $W_2^*$ such that

$$(\hat{\Phi}^*,v_n^*,v_{n,1}^*,v_{n,2}^*)\stackrel{d}{=}(\hat{\Phi},v_n,v_{n,1},v_{n,2}),$$

$$(W_\Lambda^*(C.),Z^*,W_R^*,W_1^*,W_2^*)\stackrel{d}{=}(W_\Lambda(C.),Z,W_R,W_1,W_2)$$

and for $0\leq\eta<1/2$,

(3.7) $\quad D_n:=\sup_{0\leq\theta\leq\pi/2}|\sqrt{k}(\hat{\Phi}^*(\theta)-\Phi(\theta))-(W_\Lambda^*(C_\theta)+Z^*(\theta))|=o_P(1),$

(3.8) $$\sup_{0<x,y\leq T}\frac{|v_n^*(x,y)-W_R^*(x,y)|}{(x\vee y)^\eta}=o_P(1),$$

(3.9) 
$$\sup_{0<x\leq T}\frac{|v_{n,1}^*(x)-W_1^*(x)|}{x^\eta}=o_P(1),$$
$$\sup_{0<y\leq T}\frac{|v_{n,2}^*(y)-W_2^*(y)|}{y^\eta}=o_P(1),$$

as $n \to \infty$. Henceforth, we will work on this probability space, but drop the $*$ from the notation.

PROOF OF THEOREM 2.1. By Lemma 3.1,

$$\sqrt{k}(\hat{l}_1(x,y) - l(x,y))$$
$$= \begin{cases} x\sqrt{k}\left(\hat{\Phi}\left(\frac{\pi}{2}\right) - \Phi\left(\frac{\pi}{2}\right)\right) + y \int_{\pi/4}^{\arctan y/x} \frac{1}{\sin^2 \theta} \sqrt{k}(\hat{\Phi}(\theta) - \Phi(\theta))\, d\theta, \\ \qquad \text{if } y \geq x, \\ x\sqrt{k}\left(\hat{\Phi}\left(\frac{\pi}{2}\right) - \Phi\left(\frac{\pi}{2}\right)\right) - x \int_{\arctan y/x}^{\pi/4} \frac{1}{\cos^2 \theta} \sqrt{k}(\hat{\Phi}(\theta) - \Phi(\theta))\, d\theta, \\ \qquad \text{if } y < x. \end{cases}$$

First consider the case $y \geq x$:

$$\sup_{0 < x \leq y \leq 1} \left| \frac{\sqrt{k}(\hat{l}_1(x,y) - l(x,y)) - A(x,y)}{x \vee y} \right|$$
$$= \sup_{0 < x \leq y \leq 1} \frac{1}{x \vee y} \left| x\left(\sqrt{k}\left(\hat{\Phi}\left(\frac{\pi}{2}\right) - \Phi\left(\frac{\pi}{2}\right)\right) - \left(W_\Lambda(C_{\pi/2}) + Z\left(\frac{\pi}{2}\right)\right)\right) \right.$$
$$+ y \int_{\pi/4}^{\arctan y/x} \frac{1}{\sin^2 \theta} (\sqrt{k}(\hat{\Phi}(\theta) - \Phi(\theta))$$
$$\left. - (W_\Lambda(C_\theta) + Z(\theta)))\, d\theta \right|$$
$$\leq \sup_{0 < x \leq y \leq 1} \frac{xD_n}{x \vee y} + \sup_{0 < x \leq y \leq 1} \frac{yD_n}{x \vee y} \int_{\pi/4}^{\pi/2} \frac{1}{\sin^2 \theta}\, d\theta \to 0,$$

in probability as $n \to \infty$. For the case $y < x$, the proof is similar. $\square$

Let $Q_{1n}$ and $Q_{2n}$ be the empirical quantile functions of the $\{U_i\}_{i=1}^n$ and $\{V_i\}_{i=1}^n$, respectively. Define

$$\hat{R}(x,y) = \frac{1}{k} \sum_{i=1}^n I_{\{U_i < Q_{1n}(kx/n), V_i < Q_{2n}(ky/n)\}}.$$

Note that, by (1.10),

$$\hat{l}_2(x,y) = \frac{1}{k} \sum_{i=1}^n I_{\{U_i < Q_{1n}(kx/n) \text{ or } V_i < Q_{2n}(ky/n)\}}.$$

PROOF OF THEOREM 2.2. It is easily seen that $\hat{l}_2(x,y) + \hat{R}(x,y) = (\lceil kx \rceil + \lceil ky \rceil - 2)/k \leq (\lceil kx \rceil + \lceil ky \rceil)/k$, for each $x, y \in (0,1]$, almost surely.

So we have

$$\sup_{\substack{0<x,y\leq 1 \\ x\vee y\geq 1/k}} \frac{|\sqrt{k}(\hat{l}_2(x,y) - l(x,y)) + \sqrt{k}(\hat{R}(x,y) - R(x,y))|}{(x\vee y)^\eta}$$

$$\stackrel{\text{a.s.}}{=} \sup_{\substack{0<x,y\leq 1 \\ x\vee y\geq 1/k}} \frac{|\sqrt{k}((1/k)(\lceil kx\rceil + \lceil ky\rceil - 2) - (x+y))|}{(x\vee y)^\eta}$$

$$\leq k^\eta \sup_{0<x,y\leq 1} \sqrt{k}(x+y - (\lceil kx\rceil + \lceil ky\rceil)/k)$$

$$\leq 2\sqrt{k} \cdot k^{\eta-1} = 2k^{\eta-1/2} \to 0.$$

Write $S_{jn}(x) = \frac{n}{k}Q_{jn}(\frac{k}{n}x)$, $j=1,2$. Then we have

$$\sup_{\substack{0<x,y\leq 1 \\ x\vee y\geq 1/k}} \frac{|\sqrt{k}(\hat{l}_2(x,y) - l(x,y)) + W_R(x,y) - R_1(x,y)W_1(x) - R_2(x,y)W_2(y)|}{(x\vee y)^\eta}$$

$$\stackrel{\text{a.s.}}{=} \sup_{\substack{0<x,y\leq 1 \\ x\vee y\geq 1/k}} |\sqrt{k}(\hat{R}(x,y) - R(x,y)) - W_R(x,y)$$

$$+ R_1(x,y)W_1(x) + R_2(x,y)W_2(y)|(x\vee y)^{-\eta} + o(1)$$

$$= \sup_{\substack{0<x,y\leq 1 \\ x\vee y\geq 1/k}} \frac{|\sqrt{k}(\hat{R}(x,y) - R_n(S_{1n}(x), S_{2n}(y))) - W_R(x,y)|}{(x\vee y)^\eta}$$

$$+ \sup_{\substack{0<x,y\leq 1 \\ x\vee y\geq 1/k}} \frac{|\sqrt{k}(R_n(S_{1n}(x), S_{2n}(y))) - R(S_{1n}(x), S_{2n}(y)))|}{(x\vee y)^\eta}$$

$$+ \sup_{\substack{0<x,y\leq 1 \\ x\vee y\geq 1/k}} |\sqrt{k}(R(S_{1n}(x), S_{2n}(y)) - R(x,y))$$

$$+ R_1(x,y)W_1(x,y) + R_2(x,y)W_2(y)|(x\vee y)^{-\eta} + o(1)$$

$$=: D_1 + D_2 + D_3 + o(1).$$

We will show that $D_j \to 0$ in probability, $j=1,2,3$. We have

$$D_1 = \sup_{\substack{0<x,y\leq 1 \\ x\vee y\geq 1/k}} \frac{|\sqrt{k}(T_n(S_{1n}(x), S_{2n}(y)) - R_n(S_{1n}(x), S_{2n}(y))) - W_R(x,y)|}{(x\vee y)^\eta}$$

$$\leq \sup_{\substack{0<x,y\leq 1 \\ x\vee y\geq 1/k}} |\sqrt{k}(T_n(S_{1n}(x), S_{2n}(y))$$

$$- R_n(S_{1n}(x), S_{2n}(y))) - W_R(S_{1n}(x), S_{2n}(y))|$$

$$\times (S_{1n}(x) \vee S_{2n}(y))^{-\eta} \cdot \left(\frac{S_{1n}(x) \vee S_{2n}(y)}{x \vee y}\right)^{\eta}$$

$$+ \sup_{\substack{0<x,y\leq 1 \\ x \vee y \geq 1/k}} \frac{|W_R(S_{1n}(x), S_{2n}(y)) - W_R(x,y)|}{(x \vee y)^{\eta}}$$

$$\leq \sup_{0<s,t\leq 2} \frac{|v_n(s,t) - W_R(s,t)|}{(s \vee t)^{\eta}} \cdot \sup_{\substack{0<s,t\leq k/n \\ s \vee t \geq 1/n}} \left(\frac{Q_{1n}(s) \vee Q_{2n}(t)}{s \vee t}\right)^{\eta}$$

$$+ \sup_{\substack{0<x,y\leq 1 \\ x \vee y \geq 1/k}} \frac{|W_R(S_{1n}(x), S_{2n}(y)) - W_R(x,y)|}{(x \vee y)^{\eta}}$$

$$=: D_{11} \cdot D_{12} + D_{13},$$

where the last inequality holds with arbitrarily high probability. Then $D_{11} \to 0$ in probability because of (3.8) with $T = 2$. It is well known that

$$(3.10) \qquad \sup_{s \geq 1/n} \frac{Q_{jn}(s)}{s} = O_P(1), \qquad j = 1, 2$$

(see [16], page 419). Hence, $D_{11} \cdot D_{12} \to 0$, in probability. Now consider, for each $\varepsilon > 1/k$,

$$D_{13} \leq \sup_{\substack{0<x,y\leq 1 \\ x \vee y \geq \varepsilon}} \frac{|W_R(S_{1n}(x), S_{2n}(y)) - W_R(x,y)|}{\varepsilon^{\eta}}$$

$$+ \sup_{\substack{0<x,y\leq 1 \\ 1/k \leq x \vee y \leq \varepsilon}} \frac{|W_R(S_{1n}(x), S_{2n}(y))|}{(S_{1n}(x) \vee S_{2n}(y))^{\eta}} \cdot \sup_{s,t \geq 1/n} \left(\frac{Q_{1n}(s) \vee Q_{2n}(t)}{s \vee t}\right)^{\eta}$$

$$+ \sup_{\substack{0<x,y\leq 1 \\ 1/k \leq x \vee y \leq \varepsilon}} \frac{|W_R(x,y)|}{(x \vee y)^{\eta}}$$

$$=: D_{14} + D_{15} + D_{16}.$$

By the (uniform) continuity of $W_R$ and the fact that

$$(3.11) \qquad \sup_{0<t\leq k/n} \frac{n}{k}|Q_{jn}(t) - t| \to 0, \qquad \text{a.s.}, j = 1, 2,$$

$D_{14} \to 0$ in probability a.s. for any $\varepsilon > 0$. Let $\delta > 0$; by (3.10) and Lemma 3.2, we see that, for large $n$, $P(D_{15} \geq \delta) \leq \delta$ for $\varepsilon > 0$ small enough. Again, from Lemma 3.2, we have $P(D_{16} \geq \delta) \leq \delta$. Hence, $D_{13} \to 0$ in probability and, consequently $D_1 \to 0$, in probability.

Consider $D_2$. Take $(a,b)$ with $a \vee b = u$. Then according to (2.5),

$$\frac{1}{t}C(ta, tb) = \frac{u}{ut}C\left(tu\frac{a}{u}, tu\frac{b}{u}\right)$$

$$= uR\left(\frac{a}{u}, \frac{b}{u}\right) + u^{1+\alpha}O(t^\alpha) = R(a,b) + (a \vee b)^{1+\alpha}O(t^\alpha).$$

Now with arbitrarily high probability,

$$D_2 \leq \sup_{0<x,y\leq 2} \frac{|\sqrt{k}(R_n(x,y) - R(x,y))|}{(x \vee y)^\eta} \cdot \sup_{s \vee t \geq 1/n} \left(\frac{Q_{1n}(s) \vee Q_{2n}(t)}{s \vee t}\right)^\eta.$$

We have seen before that the second term of this product is $O_P(1)$. So it suffices to show that the first term is $o(1)$:

$$\sup_{0<x,y\leq 2} \frac{|\sqrt{k}(R_n(x,y) - R(x,y))|}{(x \vee y)^\eta} = \left(\sup_{0<x,y\leq 2} \frac{\sqrt{k}(x \vee y)^{1+\alpha}}{(x \vee y)^\eta}\right) O\left(\left(\frac{k}{n}\right)^\alpha\right)$$

$$= O\left(\frac{k^{\alpha+1/2}}{n^\alpha}\right) = o(1),$$

by assumption. Hence, $D_2 \to 0$ in probability.

It remains to show that $D_3 \to 0$ in probability. By two applications of the mean-value theorem, we obtain

$$R(S_{1n}(x), S_{2n}(y)) - R(x,y)$$
$$= R(S_{1n}(x), S_{2n}(y)) - R(x, S_{2n}(y)) + R(x, S_{2n}(y)) - R(x,y)$$
$$= R_1(\theta_{1n}, S_{2n}(y))(S_{1n}(x) - x) + R_2(x, \theta_{2n})(S_{2n}(y) - y),$$

with $\theta_{1n}$ between $x$ and $S_{1n}(x)$ and $\theta_{2n}$ between $y$ and $S_{2n}(y)$. So

$$D_3 \leq \sup_{\substack{0<x,y\leq 1 \\ x \vee y \geq 1/k}} \frac{|R_1(\theta_{1n}, S_{2n}(y))\sqrt{k}(S_{1n}(x) - x) + R_1(x,y)W_1(x)|}{(x \vee y)^\eta}$$

$$+ \sup_{\substack{0<x,y\leq 1 \\ x \vee y \geq 1/k}} \frac{|R_2(x, \theta_{2n})\sqrt{k}(S_{2n}(y) - y) + R_2(x,y)W_2(y)|}{(x \vee y)^\eta}.$$

We consider only the first term on the right-hand side of this expression; the second one can be dealt with similarly. Write $z_n(x) = \sqrt{k}(S_{1n}(x) - x)$. From (3.9) with $\eta = 0$, it follows that $\sup_{0<x\leq 1} |z_n(x) + W_1(x)| \to 0$ in probability. From this, it can be shown that, for $0 \leq \eta < 1/2$,

(3.12) $$\sup_{1/k \leq x \leq 1} \frac{|z_n(x) + W_1(x)|}{x^\eta} \to 0$$

in probability (see, e.g., [7]). Now

$$\sup_{\substack{0<x,y\leq 1 \\ x \vee y \geq 1/k}} \frac{|R_1(\theta_{1n}, S_{2n}(y))z_n(x) + R_1(x,y)W_1(x)|}{(x \vee y)^\eta}$$

$$\leq \sup_{0<x,y\leq 1} R_1(\theta_{1n}, S_{2n}(y)) \cdot \sup_{1/k\leq x\leq 1} \frac{|z_n(x) + W_1(x)|}{x^\eta}$$

$$+ \sup_{0<x,y\leq 1} |R_1(x,y) - R_1(\theta_{1n}, S_{2n}(y))| \cdot \sup_{0<x\leq 1} \frac{|W_1(x)|}{x^\eta}$$

$$=: D_{31} + D_{32}.$$

Since $R_1$ is continuous on $[0,2]^2$, it is uniformly continuous and bounded. This, together with (3.12), yields $D_{31} \to 0$ in probability. The uniform continuity of $R_1$, together with (3.11) and the fact that $\sup_{0<x\leq 1}|W_1(x)|/x^\eta < \infty$ a.s., yields $D_{32} \to 0$ in probability and, consequently, $D_3 \to 0$ in probability.

Finally, we show that

$$\sup_{0<x,y<1/k} \frac{|\sqrt{k}(\hat{l}_2(x,y) - l(x,y)) + B(x,y)|}{(x \vee y)^\eta} = o_P(1).$$

Observing that $\sup_{0<x,y<1/k} \hat{l}_2(x,y) = 0$ a.s., this follows easily. □

PROOF OF THEOREM 2.3. For each $0 \leq \beta < 3$, there exist $\alpha \in [0,2)$ and $\eta \in [0, 1/2)$ such that $\beta = \alpha + 2\eta$. By Theorems 2.1 and 2.2, and $\int_0^1 \int_0^1 (x \vee y)^{-\alpha} dx\, dy < \infty$, it follows that, as $n \to \infty$,

$$\iint_{0<x,y\leq 1} \frac{k(\hat{l}_1(x,y) - \hat{l}_2(x,y))^2}{(x \vee y)^\beta} dx\, dy$$

$$= o_P(1) \iint_{0<x,y\leq 1} \frac{1}{(x \vee y)^\alpha} dx\, dy$$

$$+ \iint_{0<x,y\leq 1} \frac{(A(x,y) + B(x,y))^2}{(x \vee y)^\beta} dx\, dy$$

$$\xrightarrow{d} \iint_{0<x,y\leq 1} \frac{(A(x,y) + B(x,y))^2}{(x \vee y)^\beta} dx\, dy. \quad \Box$$

**4. Approximating the limit.** For testing purposes, we have to find the probability distribution of the limiting random variable in Theorem 2.3. This can be done by simulating the processes $A$ and $B$, but unfortunately their distributions depend on the unknown measure $\Lambda$. Therefore, we generate approximations $A_n$ and $B_n$, respectively, of the processes $A$ and $B$, not with parameter $\Lambda$, but with approximated parameter $\Lambda_n$. In this section we consider the convergence of the sequence of these approximated limiting random variables. Until further notice, we take $\{\Lambda_n\}_{n\geq 1}$ to be a sequence of *deterministic* measures.

Define

$$R_{1n}(x,y) := \tfrac{1}{2}k^{1/5}\Lambda_n([x-k^{-1/5}, x+k^{-1/5}] \times [0,y)),$$
$$R_{2n}(x,y) := \tfrac{1}{2}k^{1/5}\Lambda_n([0,x) \times [y-k^{-1/5}, y+k^{-1/5}]),$$
$$W_{R_n}(x,y) := W_{\Lambda_n}([0,x] \times [0,y]),$$
$$W_{1n}(x) := W_{\Lambda_n}([0,x] \times [0,\infty]),$$
$$W_{2n}(y) := W_{\Lambda_n}([0,\infty] \times [0,y]),$$

and the process $B_n$ by

$$B_n(x,y) := W_{R_n}(x,y) - R_{1n}(x,y)W_{1n}(x) - R_{2n}(x,y)W_{2n}(y).$$

Based on the definition of $Z$ in (2.1) and the homogeneity property of $\lambda$ [i.e., $\lambda(tx,ty) = \tfrac{1}{t}\lambda(x,y)$], we define the approximating process $Z_n$ by

$$(4.1) \quad Z_n(\theta) = \begin{cases} \lambda_n(1, \tan\theta) \tan\theta \int_0^{1/\tan\theta} \frac{W_{1n}(x)}{x} dx \\ \quad - \lambda_n(1, \tan\theta) \int_0^1 \frac{W_{2n}(x)}{x} dx \\ \quad - W_{2n}(1) \int_{1/\tan\theta}^{\infty} \lambda_n(x,1) dx, \qquad \theta \in [0, \pi/4], \\ \lambda_n(1/\tan\theta, 1) \int_0^1 \frac{W_{1n}(x)}{x} dx \\ \quad - \lambda_n(1/\tan\theta, 1) \frac{1}{\tan\theta} \int_0^{\tan\theta} \frac{W_{2n}(x)}{x} dx \\ \quad - W_{2n}(1) \int_1^{\infty} \lambda_n(x,1) dx - W_{1n}(1) \int_1^{\tan\theta} \lambda_n(1,y) dy, \\ \qquad\qquad\qquad\qquad\qquad\qquad\qquad \theta \in (\pi/4, \pi/2), \\ -W_{2n}(1) \int_1^{\infty} \lambda_n(x,1) dx - W_{1n}(1) \int_1^{\infty} \lambda_n(1,y) dy, \\ \qquad\qquad\qquad\qquad\qquad\qquad\qquad \theta = \pi/2, \end{cases}$$

where $\lambda_n$ is the approximation of $\lambda$ defined by

$$\lambda_n(1,y) := \tfrac{1}{4}k^{1/3}\Lambda_n([1-k^{-1/6}, 1+k^{-1/6}] \times [y-k^{-1/6}, y+k^{-1/6}]),$$
$$y > 0,$$
$$\lambda_n(x,1) := \tfrac{1}{4}k^{1/3}\Lambda_n([x-k^{-1/6}, x+k^{-1/6}] \times [1-k^{-1/6}, 1+k^{-1/6}]),$$
$$x > 0.$$

Finally, define the process $A_n$ by

$$A_n(x,y) := \begin{cases} x\left(W_{\Lambda_n}(C_{\pi/2}) + Z_n\left(\frac{\pi}{2}\right)\right) \\ \quad + y \int_{\pi/4}^{\arctan y/x} \frac{1}{\sin^2\theta}(W_{\Lambda_n}(C_\theta) + Z_n(\theta))\,d\theta, & \text{if } y \geq x, \\ x\left(W_{\Lambda_n}(C_{\pi/2}) + Z_n\left(\frac{\pi}{2}\right)\right) \\ \quad - x \int_{\arctan y/x}^{\pi/4} \frac{1}{\cos^2\theta}(W_{\Lambda_n}(C_\theta) + Z_n(\theta))\,d\theta, & \text{if } y < x. \end{cases}$$

First we consider the weak convergence of the weighted approximating processes. We write $D_2 := D[0,1]^2$ for the generalization of $D[0,1]$ to dimension 2, and $\mathcal{L}_d$ for the Borel $\sigma$-algebra on $(D_2, d)$, where $d$ is the metric on $D_2$ defined in [11].

PROPOSITION 4.1. *Let $\Lambda$ be as in Theorem 2.3. Suppose that $\{\Lambda_n\}_{n\geq 1}$ is a sequence of measures on $[0,\infty]^2 \setminus \{(\infty,\infty)\}$ satisfying that, for each $x, y \geq 0$,*

(4.2) $\quad \Lambda_n([0,x] \times [0,\infty]) = [kx]/k, \qquad \Lambda_n([0,\infty] \times [0,y]) = [ky]/k,$

(4.3) $\quad \sup_{0<x,y\leq 1} |\Lambda_n([0,x]\times[0,y]) - \Lambda([0,x]\times[0,y])| \to 0$

*as $n \to \infty$. Further, suppose that*

(4.4) $\quad \sup_{0<x\leq 1} |\lambda_n(x,1) - \lambda(x,1)| \to 0, \qquad \sup_{0<y\leq 1} |\lambda_n(1,y) - \lambda(1,y)| \to 0,$

(4.5) $\quad \sup_{0<x,y\leq 1} |R_{jn}(x,y) - R_j(x,y)| \to 0, \qquad j = 1,2,$

*as $n \to \infty$. Then for each $0 \leq \eta < 1/2$,*

$$\left\{\frac{A_n(x,y) + B_n(x,y)}{(x \vee y)^\eta}, (x,y) \in [0,1]^2\right\} \to \left\{\frac{A(x,y) + B(x,y)}{(x \vee y)^\eta}, (x,y) \in [0,1]^2\right\},$$

*weakly in $D_2$.*

Before proving this proposition, we present three corollaries. The last one is the main result of this section.

COROLLARY 4.1. *Under the conditions of Proposition 4.1, for each $0 \leq \beta < 3$,*

(4.6)
$$\iint_{0<x,y\leq 1} \frac{(A_n(x,y) + B_n(x,y))^2}{(x \vee y)^\beta}\,dx\,dy$$
$$\xrightarrow{d} \iint_{0<x,y\leq 1} \frac{(A(x,y) + B(x,y))^2}{(x \vee y)^\beta}\,dx\,dy$$

*as* $n \to \infty$.

Let $Q_{\Lambda_n}$ be the quantile function of the random variable on the left-hand side of (4.6) and $Q_\Lambda$ the quantile function of the random variable on the right-hand side of (4.6).

COROLLARY 4.2. *Under the conditions of Proposition* 4.1, *for each* $0 \leq \beta < 3$ *and for each continuity point* $1 - \alpha$ $(0 < \alpha < 1)$ *of* $Q_\Lambda$,

$$\lim_{n \to \infty} Q_{\Lambda_n}(1 - \alpha) = Q_\Lambda(1 - \alpha).$$

Next, with abuse of notation, we estimate $\Lambda_n$ from the data, so it becomes random. In [8], $\Lambda_n$ is defined as

(4.7)
$$\Lambda_n(A) := \frac{1}{k} \sum_{i=1}^n I_{kA/n}\left(\frac{1}{n}\sum_{j=1}^n I_{(-\infty, U_i]}(U_j), \frac{1}{n}\sum_{j=1}^n I_{(-\infty, V_i]}(V_j)\right)$$
$$= \frac{1}{k}\sum_{i=1}^n I_{kA}(n+1-R_i^X, n+1-R_i^Y),$$

where $U_i := 1 - F_1(X_i)$, $V_i := 1 - F_2(Y_i)$, for $i = 1, 2, \ldots, n$. Note that, for $x, y > 0$,

$$\Lambda_n([0,x) \times [0,y)) = \frac{1}{k}\sum_{i=1}^n I_{\{U_i < Q_{1n}(kx/n), V_i < Q_{2n}(ky/n)\}}.$$

So $\Lambda_n([0, x) \times [0, \infty]) = (\lceil kx \rceil - 1)/k \leq [kx]/k = \Lambda_n([0, x] \times [0, \infty])$ a.s. and $\Lambda_n([0, \infty] \times [0, y)) = (\lceil ky \rceil - 1)/k \leq [ky]/k = \Lambda_n([0, \infty] \times [0, y])$ a.s.

The final and main corollary deals with the *random* measures $\Lambda_n$, where the functions derived from $\Lambda_n$, like $\lambda_n$, are defined as before. In particular, we define $Q_{\Lambda_n}$ as the quantile function of the random variable on the left-hand side of (4.6), *conditional* on $\Lambda_n$, so it is also random.

COROLLARY 4.3. *Let* $\Lambda_n$ *be as in* (4.7). *Under the conditions of Theorem* 2.3, *we have, for each* $0 \leq \beta < 3$ *and each continuity point* $1 - \alpha$ $(0 < \alpha < 1)$ *of* $Q_\Lambda$, *that*

$$Q_{\Lambda_n}(1 - \alpha) \xrightarrow{P} Q_\Lambda(1 - \alpha) \qquad \text{as } n \to \infty.$$

For testing purposes, Corollary 4.3 shows that simulation of the limiting random variable in Theorem 2.3 with $\Lambda$ replaced by the estimated $\Lambda_n$ is asymptotically correct.

Now we turn to the proofs. In order to prove Proposition 4.1, by Prohorov's theorem, it is necessary and sufficient to prove the following:

(i) The finite-dimensional distributions of $\{(A_n(x,y) + B_n(x,y))/(x \vee y)^\eta, (x,y) \in [0,1]^2\}_{n \geq 1}$ converge to those of $\{(A(x,y) + B(x,y))/(x \vee y)^\eta, (x,y) \in [0,1]^2\}$;

(ii) $\{(A_n(x,y) + B_n(x,y))/(x \vee y)^\eta, (x,y) \in [0,1]^2\}_{n \geq 1}$ is relatively compact.

For the relative compactness, we need several lemmas. These lemmas and their proofs can be found in a separate Appendix, posted at http://center.uvt.nl/staff/einmahl/AppEdHL.pdf, or in [10], pages 81–87. These lemmas lead to the following results: Under the conditions of Proposition 4.1, for each $0 \leq \eta < 1/2$,

$$\left\{\frac{B_n(x,y)}{(x \vee y)^\eta}, (x,y) \in [0,1]^2\right\}_{n \geq 1}$$

is relatively compact, and for each $0 \leq \eta < 1$,

$$\left\{\frac{A_n(x,y)}{(x \vee y)^\eta}, (x,y) \in [0,1]^2\right\}_{n \geq 1}$$

is relatively compact.

PROOF OF PROPOSITION 4.1. By these results,

(4.8) $$\left\{\frac{A_n(x,y) + B_n(x,y)}{(x \vee y)^\eta}, (x,y) \in [0,1]^2\right\}_{n \geq 1}$$

is relatively compact. It is easy to check that the finite-dimensional distributions of our estimated processes in (4.8) converge to those of the limiting process, which completes the proof. □

PROOF OF COROLLARY 4.1. After applying a Skorohod construction to the weak convergence statement of Proposition 4.1, the proof is similar to that of Theorem 2.3. □

PROOF OF COROLLARY 4.2. Proposition 4.1 implies the weak convergence of the distribution function of the left-hand side of (4.6) to the distribution function of the right-hand side of (4.6). This property carries over to the inverse functions $Q_{\Lambda_n}$ and $Q_\Lambda$. □

PROOF OF COROLLARY 4.3. From another Skorohod construction, we obtain an a.s. version of the statement of Theorem 2.2; without changing the notation, we now work with this construction. Since, for $0 < x, y \leq 1$,

$$\Lambda([0,x] \times [0,y]) = x + y - l(x,y),$$
$$\Lambda_n([0,x] \times [0,y]) = \lceil kx \rceil/k + \lceil ky \rceil/k - \hat{l}_2(x,y) - \delta_n(x,y)/k$$

[$\delta_n(x,y)$ takes values in $\{0,1,2\}$], it follows that, for each $\varepsilon > 0$,

$$(4.9) \quad \sup_{0<x,y\leq 1} k^{1/2-\varepsilon}|\Lambda_n([0,x]\times[0,y]) - \Lambda([0,x]\times[0,y])| \to 0 \quad \text{a.s.,}$$

as $n \to \infty$.

We now show that (4.2), (4.3), (4.4) and (4.5) hold a.s. We already saw, below (4.7), that (4.2) holds a.s. and the a.s. version of (4.3) follows immediately from (4.9).

By (4.9) and (4.2), it easily follows that

$$(4.10) \quad \sup_{E \in \mathcal{E}} k^{1/2-\varepsilon}|\Lambda_n(E) - \Lambda(E)| \to 0 \quad \text{a.s.,}$$

as $n \to \infty$, where $\mathcal{E} := \{E | E = [x_1, x_2] \times [y_1, y_2], 0 < x_1 \leq x_2 \leq 2, 0 < y_1 \leq y_2 \leq 2\}$. Let $E_n(x) = [x - k^{-1/6}, x + k^{-1/6}] \times [1 - k^{-1/6}, 1 + k^{-1/6}]$. Then, setting $\lambda(u,v) = 0$ if $u < 0$,

$$\sup_{0<x\leq 1} |\lambda_n(x,1) - \lambda(x,1)|$$

$$\leq \sup_{0<x\leq 1} \tfrac{1}{4}k^{1/3}|\Lambda_n(E_n(x)) - \Lambda(E_n(x))| + \sup_{0<x\leq 1} |\tfrac{1}{4}k^{1/3}\Lambda(E_n(x)) - \lambda(x,1)|$$

$$\leq \sup_{0<x\leq 1} \tfrac{1}{4}k^{1/3}|\Lambda_n(E_n(x)) - \Lambda(E_n(x))|$$

$$+ \sup_{0<x\leq 1} \sup_{(u,v)\in E_n(x)} |\lambda(u,v) - \lambda(x,1)| \to 0 \quad \text{a.s.,}$$

as $n \to \infty$, by (4.10) and the uniform continuity of $\lambda$ on $[-1,2] \times [\tfrac{1}{2}, 2]$ [which follows from $\lambda(0,1) = 0$]. The proofs of $\sup_{0<y\leq 1}|\lambda_n(1,y) - \lambda(1,y)| \to 0$ a.s. and $\sup_{0<x,y\leq 1}|R_{jn}(x,y) - R_j(x,y)| \to 0$, $j = 1, 2$, a.s. are similar. Hence, (4.4) and (4.5) hold a.s.

According to Corollary 4.2 we have $Q_{\Lambda_n}(1-\alpha) \to Q_\Lambda(1-\alpha)$ a.s., as $n \to \infty$, hence, also in probability. $\square$

**5. Simulation study and real data application.** In this section we present a simulation study, making use of the results of Section 4. We will consider two distributions satisfying the domain of attraction condition and one that fails to satisfy it. At the end of the section we will apply our procedure to financial data.

Theoretically, we can choose any $\beta \in [0,3)$ in the test statistic in (1.12). We investigate the influence of $\beta$ on the testing procedure by sampling from the bivariate Cauchy distribution. We choose $\beta$ to be 0, 1 or 2.

Consider the bivariate Cauchy distribution restricted to the first quadrant, with density

$$f(x,y) = \frac{2}{\pi(1 + x^2 + y^2)^{3/2}}, \qquad x, y > 0.$$

TESTING THE BIVARIATE EVT CONDITION 25TABLE 1
*Quantiles of the limiting r.v. for the bivariate Cauchy distribution*

| $\beta$ | $p$ | | | | | | | |
|---|---|---|---|---|---|---|---|---|
| | 0.10 | 0.25 | 0.50 | 0.75 | 0.90 | 0.95 | 0.975 | 0.99 |
| 0 | 0.018 | 0.025 | 0.038 | 0.065 | 0.106 | 0.142 | 0.177 | 0.227 |
| 1 | 0.030 | 0.041 | 0.062 | 0.103 | 0.168 | 0.222 | 0.278 | 0.356 |
| 2 | 0.074 | 0.099 | 0.144 | 0.224 | 0.347 | 0.447 | 0.554 | 0.699 |

It readily follows that

$$\Lambda([0,x] \times [0,y]) = x + y - \sqrt{x^2 + y^2},$$

$$\lambda(x,y) = \frac{xy}{(x^2+y^2)^{3/2}}, \qquad x, y > 0.$$

This distribution satisfies the conditions of Theorem 2.3; in particular, (2.5) holds with $\alpha = 2$ (see [8], pages 1409–1410). First we present in Table 1 the quantiles of the limiting random variable

$$\iint_{0 < x,y \leq 1} \frac{(A(x,y) + B(x,y))^2}{(x \vee y)^\beta} \, dx \, dy.$$

We used 1,000,000 replications. With high probability, these quantiles are accurate up to 0.01.

Now for sample size $n = 2000$, we simulated 2000 times the test statistic $kL_n$, for various values of $k$. Using the 0.95th quantiles above, we find the simulated type-I error probabilities; see Table 2. In this table also the empirical median and the empirical 0.95th quantile of the test statistics are shown. In the ideal situation the number of rejections is a binomial r.v. with parameters 2000 and 0.05. So the simulated type-I errors in Table 2 are remarkably close to 0.05. Only for $k = 400$ does bias seem to set in. Also, the empirical median and 0.95th quantile of the test statistics are very close to those of the limiting r.v. listed in Table 1. Generally speaking, the influence of $\beta$ on the quality of the results is very small for the Cauchy distribution. From Table 2, we feel that $\beta = 2$ works slightly better than the others. Because of this and because we want to put additional emphasis on the extreme observations, from now on we take $\beta = 2$.

In practice, for a given dataset, we first calculate the test statistic $kL_n$; then we estimate the measure $\Lambda_n$ and simulate the 0.95th quantile of the estimated limiting r.v. using the approximation of Section 4. Finally, if the test statistic is not smaller than this 0.95th quantile, we reject the null hypothesis (1.8).

First we consider again the bivariate Cauchy distribution and take two samples of size $n = 2000$. The results are presented in Figure 1. Note that the behavior of the test statistic and the estimated 0.95th quantile fluctuate with the sample fraction $k$, but that, for all $k$ in the figure, the value is— correctly—far below the estimated 0.95th quantile of the limiting random variable.

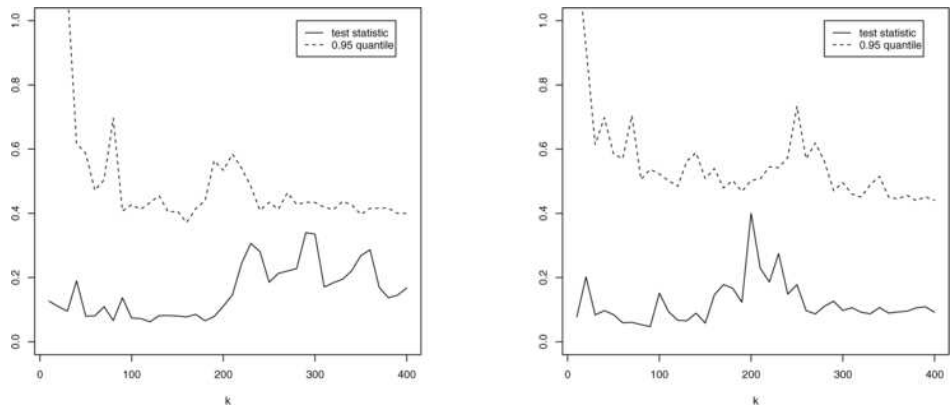

Fig. 1.  *Cauchy distribution: test statistic and* 0.95*th quantile for two samples.*

Next we generate 2000 independent pairs $(U, 1-V)$, where $(U,V)$ has a Gumbel copula as distribution function, that is, the d.f. is given by

$$C(u,v) = \exp(-[(-\log u)^\theta + (-\log v)^\theta]^{1/\theta}), \qquad \theta \geq 1;$$

we take $\theta = 10$. It is easily checked that, for the d.f. of $(U, 1-V)$, (1.1) holds and that we have asymptotic independence; see Remark 2.2. Since our results do not apply for the case of asymptotic independence, we only present the test statistic itself (Figure 2, left panel). We see that, for $k$ up to 200, the test statistic is very close to 0 (which strongly supports $H_0$) and that bias sets in

Table 2
*Simulated type-I error, median and* 0.95*th quantile of the test statistics for the Cauchy d.f.;* $n = 2000$, $\alpha = 0.05$

| $\beta$ | $k$ | 20 | 40 | 60 | 80 | 100 | 125 | 150 | 175 | 200 | 300 | 350 | 400 |
|---|---|---|---|---|---|---|---|---|---|---|---|---|---|
| 0 | $\hat{\alpha}$ | 0.041 | 0.045 | 0.047 | 0.044 | 0.038 | 0.047 | 0.047 | 0.034 | 0.035 | 0.049 | 0.048 | 0.060 |
|   | Q(0.95) | 0.134 | 0.135 | 0.139 | 0.132 | 0.129 | 0.139 | 0.139 | 0.127 | 0.125 | 0.141 | 0.140 | 0.153 |
|   | Q(0.50) | 0.036 | 0.038 | 0.036 | 0.036 | 0.036 | 0.036 | 0.036 | 0.036 | 0.036 | 0.040 | 0.040 | 0.047 |
| 1 | $\hat{\alpha}$ | 0.041 | 0.047 | 0.047 | 0.045 | 0.039 | 0.050 | 0.044 | 0.036 | 0.034 | 0.054 | 0.046 | 0.061 |
|   | Q(0.95) | 0.208 | 0.213 | 0.216 | 0.210 | 0.210 | 0.220 | 0.216 | 0.203 | 0.203 | 0.226 | 0.216 | 0.236 |
|   | Q(0.50) | 0.059 | 0.061 | 0.059 | 0.059 | 0.059 | 0.060 | 0.059 | 0.058 | 0.058 | 0.065 | 0.064 | 0.076 |
| 2 | $\hat{\alpha}$ | 0.047 | 0.042 | 0.049 | 0.048 | 0.044 | 0.047 | 0.044 | 0.044 | 0.042 | 0.053 | 0.050 | 0.068 |
|   | Q(0.95) | 0.434 | 0.423 | 0.444 | 0.442 | 0.430 | 0.437 | 0.431 | 0.431 | 0.416 | 0.463 | 0.446 | 0.503 |
|   | Q(0.50) | 0.133 | 0.138 | 0.137 | 0.135 | 0.137 | 0.141 | 0.138 | 0.143 | 0.143 | 0.156 | 0.156 | 0.195 |

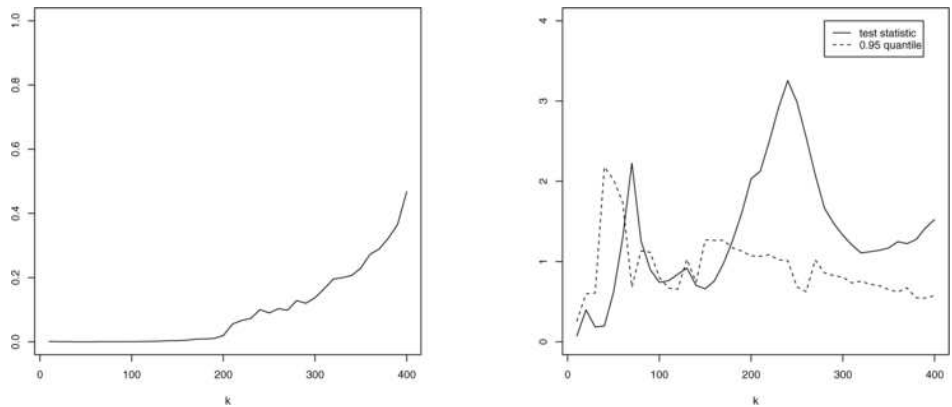

Fig. 2. *Test statistic transformed Gumbel copula (*left*) and alternative distribution (*right*).*

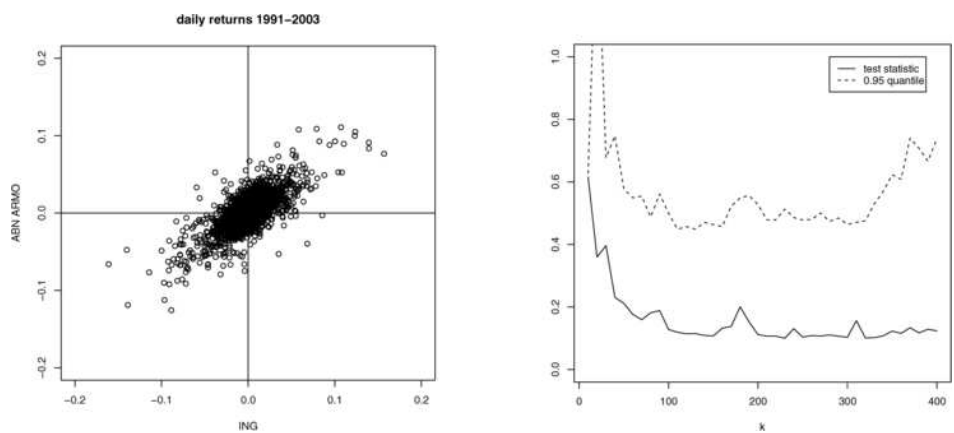

Fig. 3. *Daily equity returns of two Dutch banks (*left*) and test statistic and* 0.95*th quantile (*right*).*

for larger values of $k$. We also consider 2000 observations from a d.f. (which is also a copula), which does not satisfy condition (1.8). The distribution is an adaptation of a distribution in [15]: take a density of $3/2$ on the following rectangles: $[2^{-(2m+1)}, 2^{-(2m)}] \times [2^{-(2r+1)}, 2^{-(2r)}]$, for $m = 0, 1, 2, \ldots$ and $r = 0, 1, 2, \ldots$; in this way a probability mass of $2/3$ is assigned. The remaining $1/3$ is assigned by taking the uniform distribution on the line segments from $(2^{-(2m+2)}, 2^{-(2m+2)})$ to $(2^{-(2m+1)}, 2^{-(2m+1)})$, $m = 0, 1, 2, \ldots$, such that the mass of the $m$th segment is equal to $2^{-(2m+2)}$. In Figure 2 (right panel), we see, for varying $k$, the test statistic and simulated 0.95th quantile of the sample of size $n = 2000$ from this distribution. Again, the test statistic and

28 J. H. J. EINMAHL, L. DE HAAN AND D. LI

the estimated 0.95th quantile fluctuate with $k$, but from a certain $k$ on (and for most values of $k$), the null hypothesis is clearly rejected.

Finally, we apply the test to real data, similarly as we just did for the simulated data sets in Figures 1 and 2. The data are 3283 daily logarithmic equity returns over the period 1991–2003 for two Dutch banks, ING and ABN AMRO bank. The bivariate, heavy-tailed data are shown in Figure 3 on the left; on the right we see again the test statistic and the simulated 0.95th quantile. Since the test statistic is everywhere clearly below the quantile, we cannot reject the null hypothesis. This is a satisfactory result, because it allows us to analyze these data further, using the statistical theory of extremes.

**Acknowledgment.** We are grateful to two referees for a careful reading of the manuscript and several useful comments.

J. H. J. EINMAHL
DEPARTMENT OF ECONOMETRICS AND CENTER
TILBURG UNIVERSITY
P.O. BOX 90153
5000 LE TILBURG
THE NETHERLANDS
E-MAIL: j.h.j.einmahl@uvt.nl

L. DE HAAN
ECONOMETRIC INSTITUTE
ERASMUS UNIVERSITY
P.O. BOX 1738
3000 DR ROTTERDAM
THE NETHERLANDS
E-MAIL: ldehaan@few.eur.nl

D. LI
IMSV
UNIVERSITY OF BERN
SIDLERSTRASSE 5
CH-3012 BERN
SWITZERLAND
E-MAIL: deyuan.li@stat.unibe.ch